\documentclass[12pt,a4paper]{article}
\usepackage{amssymb,amsmath}
\newcommand{\cqfd}{ \vskip-5mm \mbox{}
\nolinebreak \hfill $\Box$ \vskip 3mm}

\def\T{{\mathbb T}}
\def\R{\mathbb R}

\title{Examples of quasi-hyperbolic dynamical systems with slow decay of correlations} %titre dans l'autre langue

\author{St\'ephane Le Borgne}
\date{Universit\'e de Rennes, Campus de Beaulieu, 35042 Rennes Cedex}
\begin{document}

    \maketitle

{\bf Abstract} :\  We give examples  of quasi-hyperbolic dynamical systems with the following properties : polynomial decay of correlations, convergence in law toward a non gaussian law of the ergodic sums (divided by $n^{3/4}$) associated to non degenerated regular functions.

{\bf Keywords} :\  quasi-hyperbolicity, partial hyperbolicity, 
limit theorem.

\vskip 1cm

In the litterature dealing with stochastic properties
 of   quasi-hyperbolic transformations of manifolds
 one establish most often exponential decay of correlations and
 central limit theorem. For a few years several examples of dynamical systems were found for which the correlations and the behaviour in law are different (\cite{4}, \cite{10}, \cite{11}, \cite{13}, \cite{15}). Here we present a new family of examples (of systems that can be invertible) constructed in a very simple way as skew products. 

Let $(X,A,\mu)$ be an Anosov system or an Axiom-A system and  $(Y,\phi_t,\nu)$ an Anosov flow both with exponential decay of correlations. Let us fix a regular function $f$ from  $X$ to $\mathbb R$.
On the  product $X\times Y$ d\'efine the transformation
 \begin{eqnarray*} T :  X\times Y&\longrightarrow
 &X\times Y\\
 (x,y)& \longmapsto& (A x, \phi_{f(x)} y).
 \end{eqnarray*}
The hyperbolic transformation  $A$ acts on the first coordinate ; the flow acts  on the second one during a time which depends on the first coordinate via $f$. The  product measure $\mu\otimes\nu$  is $T$-invariant :  $(X\times Y, T, \mu\otimes\nu)$ is a dynamical system.

Let $S_n f(x)$ denote $\displaystyle\mathop{\sum}_{k=0}^{n-1} f(A^k x)$ 
($S_0=0$) the ergodic sums associated to the action of $A$ on the function $f$. A simple computation 
shows that, for all $n\geq 0$, the expression of $T^n$  is:
$$T^n(x,y)=(A^n x, g_{S_n f(x)} y)$$
As $f$ is centered the ergodic $S_n f$ sums grows almost surely less rapidly than $n^{1/2+\varepsilon}$ (for all $\varepsilon>0$). So transformation $T$ has a null Lyapounoff exponent in the tangent spaces of the submanifold $\{x\}\times Y$.

\section{ Decay of correlation}

Let $\sigma(f)$ denote the non-negative quantity defined by
$$
\sigma^2(f)=\int_{X}f^2d\mu+2\sum_{k\geq 1}<f,A^kf>.
$$ 

{\bf Theorem} :\it   If the function $f$ is aperiodic and if $\varphi : X\times Y \rightarrow \mathbb R$ is a centered  h\"older-continuous function such that the function $y\longmapsto \int_{X}
\varphi(a,y) d\mu(a)$ is not a coboundary  in the  sense of the flow the the number
$$
\Sigma^2(\varphi)=\int_{\mathbb  R}\int_{Y}
\left(\int_{X}\varphi(a,g_b y) d\mu(a)\right) \left(\int_{X} \varphi (c,y) d\mu(c)\right) d\nu(y) \ db
$$
is positive and the $<\varphi\circ T^k,\varphi>$ is equivalent to 
$$
\frac{1}{\sqrt{2\pi k} \sigma(f)}\sum^2(\varphi).
$$
\rm

{\it Proof} : We're interested in
\begin{eqnarray*} <\varphi\circ T^k,\varphi>&=&\mathbb E (\varphi(T^k.,g_{S_k f(.)}.)\varphi (.,.))\\
& =& \int_{X\times Y} \varphi (T^k x, g_{S_{k-1}f(x)}y) \varphi (x,y) d\mu(x) d\nu(y).
\end{eqnarray*}

Let us consider the  function $F$ :
\begin{eqnarray*} X\times \mathbb R\times
X&\longrightarrow&
\mathbb R\\
(a,b,c)&\longmapsto&\int_{Y} \varphi (a, g_by) \varphi (c,y) d\nu(y)-(\int_{Y} \varphi (a, y) d\nu(y))
(\int_{Y} \varphi (c,y)d\nu(y))
\end{eqnarray*}

The correlation $<\varphi\circ T^k,\varphi>$ is
$$\int_{X} F(T^{k}x, S_{k-1} f(x),x) d\mu(x)+\int_{X}
(\int_{Y} \varphi (T^k x,y) d\nu(y)) (\int_{Y}\varphi (x,y) d\nu(y)) d\mu(x).$$
The action of $A$ is exponentially mixing thus the second element of this sum tends exponentially fast toward 0.

The action of the  flow is exponentially mixing. The  function $F$ decrease exponentially fast when $b$ grows.
The local limit theorem of Guivarc'h and Hardy assure that the
integral 
$$
\int_{X} F(T^k x,S_{k-1} f(x),x) d\mu(x)
$$ is
equivalent to, 
\begin{eqnarray*}
&{ }&   \frac{1}{\sqrt{2k\pi} \sigma(f)} \int_{X\times\mathbb R\times X} F(a,b,c) d\mu(a)\ db\ d\mu(c)\\
&{ }& \ \ \ =\frac{1}{\sqrt{2\pi k} \sigma(f)} \int_{\mathbb R} \big(\int_{X}
\int_{X}\int_{Y} \varphi (a,g_by)\varphi(c,y) d\nu(y)\ d\mu(a)\ d\mu(c) \\
&{ }& \ \ \ \ \ \ -\underbrace{\int_{X}
\int_{X}\int_{Y} \varphi (a,y) d\nu(y)\int_{Y}
\varphi (c,y) d\nu(y)\ d\mu(a)\ d\mu(c)}_0 \big) db\\
&{ }& \ \ \ =\frac{1}{\sqrt{2\pi k} \sigma(f)} \underbrace{\int_{\mathbb  R}\int_{Y}
\left(\int_{X}\varphi(a,g_b y) d\mu(a)\right) \left(\int_{X} \varphi (c,y) d\mu(c)\right) d\nu(y) \ db.}_{\sum^2(\varphi).}\\
\end{eqnarray*}

The quantity $\sum^2(\varphi)$ is zero if and only if the function $y\longmapsto \int_{X}
\varphi(a,y) d\mu(a)$ is a coboundary  in the sense of the flow (see \cite{3} for more details) that is if and only if there exists a function $\varphi : Y \rightarrow \mathbb R$ such that for almost every $y\in Y$ one has
$$\int_{X} \varphi(a,y)
d\mu(a)=\displaystyle\mathop{lim}_{t\rightarrow 0} \frac{\psi (g+y)-\varphi (y)}{t}.$$
\cqfd

\section{Variances of the ergodic sums}

{\bf Corollary} :\it  If $\sum^2(\varphi)$ is not zero the variance of ${\sum}_{k=0}^{n-1}\varphi\circ T^k)$ is equivalent to 
$$
\frac{8}{3}
\frac{\sum^2(\varphi)}{\sqrt{2\pi}\sigma(f)}n^{3/2}.
$$
\rm 

{\it  Proof} : This is a direct computation using the theorem.\cqfd

{\bf Remark } : An analogous construction is possible in continuous time. Let $(Z,\psi_t,\rho)$ be another Anosov flow. Define 
 \begin{eqnarray*} \chi_t :  Z\times Y&\longrightarrow
 &Z\times Y\\
 (z,y)& \longmapsto& (\chi_t z, \phi_{\int_0^tf(\psi_s z)ds} y).
 \end{eqnarray*}
One obtain a quasi-hyperbolic flow with decay of  correlations in $t^{-1/2}$. To show this one has to use the results of Waddington \cite{14}.

\section{Convergence in law}

Once we know the behaviour of the variances, we can ask the question of the convergence in law. We treat this problem on an example. Many technical difficulties  arise when we try to generalise our result.

Our example is the following one. Let us consider

- the group $G=PSL (2,\mathbb R)$,

- $\Gamma$ a cocompact lattice in $G$,

- the one-parameter group $g_t:= \left(\begin{array}{cc}e^{t/2}&0\\0&e^{-t/2}
                                    \end{array}\right)$,

- the two dimensional torus $2,\T^2$,

- a centered function $f$ from $\T^2$ to $\mathbb R$,

 - the matrix $A:=
 \left(\begin{array}{cc}2&1\\1&1\end{array}\right)$.

In 1988 Rudolph \cite{12} has introduced the application

 \begin{eqnarray*} T : \T^2 \times G/\Gamma&\longrightarrow
 &\T^2\times G/\Gamma\\
 (x,y)& \longmapsto& (A x, g_{f(x)} y).
 \end{eqnarray*}

An example of concrete function  for which the system has good properties is the function :
\begin{eqnarray*}f :\qquad \T^2&\longrightarrow&\mathbb R\\
(x_1,x_2)&\longmapsto&sin (2\pi x_1).
 \end{eqnarray*}

The transformation $T$ has two non zero Lyapounoff exponents, $\lambda=\frac{3+\sqrt 5}
{2}$ and $\frac{1}{\lambda}$, associated to two directions one expanded the other contracted uniformly by  $T$ and a zero exponent in the tangent spaces to the submanifolds $\{x\}\times G/\Gamma$. Let $\mu$ and $\nu$ be the Haar measure and the measure coming from the Haar measure respectively.

The system $(\T^2\times G/\Gamma, T, \mu\otimes\nu)$ is a regular version of the
 $T,T^{-1}$-transformation. Both are  examples of $K$-systems that do not have the Benoulli property (\cite{7}, \cite{12}) . 

In $\mathbb R^2$ let $\tilde{x_0}$ (resp.$\tilde{x_{-1}}$   denote the intersection point of the contracted right line passing through the point $(1,0)$ (resp. $(1,-1)$) and the expanded right line passing through the point $(1,1)$. These are a homoclinic points : in the torus $T^k\tilde{x_0}$ tend and $T^k\tilde{x_{-1}}$ to zero when $k$ tends toward $+\infty$ and toward $-\infty$.

{\bf Theorem} : \it  Let  $f$ be a h\"older continuous function such that  
$$
\sum_{-\infty}^{\infty}(f(T^k\tilde{x_0})-f(T^k\tilde{x_{-1}}))\neq0.
$$
Let
$\varphi : \T^2 \times G/\Gamma\longrightarrow \mathbb R $ be a h\"older continuous function such that $\sum^2 (\varphi)$ is positive. There exist three brownian motion (non necessarily
reduced) $W, W_+,W_-$ such that, if $L_t(x)$ denotes the local time of $W$ in $x$, one has :
$$
\frac{1}{n^{3/4}} \displaystyle\mathop{\sum}_{k=0}^{n-1}
\varphi\circ T^k\displaystyle\mathop{\longrightarrow}^{\mathcal L} \int^{+\infty}_0
L_1(x) dW_+(x)+\int^{+\infty}_0 L_1(-x)dW_-(x)
$$
\rm

{\bf Remark } : One easily verify that if $f$ is a coboundary the condition on $f$ isn't satisfied. For regular function and for  the automorphism $A$, to be aperiodic simply means not to be a coboundary. The condition of the theorem  is {\sl a priori} stronger than aperiodicity. Is it really stronger ? Another question : for which systems and which functions the notions of periodicity and coboundaricity do coincide ?

One always has $\Sigma^2(\varphi-\int_{\T^ 2} \varphi(a,.) d\mu(a))=0$. The study of the convergence
in law of $\sum \varphi \circ T^k$ is thus reduced to the one of $\sum(\int_{\T^2} \varphi
(a,.) d\mu(a))\circ T^k$. So one can restrict ourselves to the case where $\varphi$ depends uniquely on the second
coordinate. That what we're doing from now on.

We'll use the following result of
Kesten and Spitzer.

{\bf Theorem} : (Kesten, Spitzer)
\it  Let $X_i$ an i.i.d sequence of centered random variables  with values in $\mathbb Z$. Let
 $(\xi_i)$ an i.i.d sequence of centered random variables  independent of the $X_i$ variables. There exist three brownian motion (non necessarily
reduced) $W, W_+,W_-$ such that, if $L_t(x)$ denotes the local time of $W$ in $x$, one has :
$$
\frac{1}{n^{3/4}} \displaystyle\mathop{\sum}_{k=0}^{n-1}
 \xi_{X_i+\ldots+X_k}\displaystyle\mathop{\longrightarrow}^{\mathcal L} \int^{+\infty}_0
 L_1(x) dW_+(x)+\int^{+\infty}_0 L_1(-x)dW_-(x)
$$

\rm

The theorem is a direct consequence of the two following propositions and of the result of Kesten and Spitzer.

Define
$$
N(n,p)=\sharp\{k\in\{0,\ldots,n-1\}|/S_k f\in [p,p+1[\}.
$$

{\bf  Proposition} : \it  Under the hypotheses of the theorem, one has
$$
\frac{1}{n^{3/4}} \sum_{p\in \mathbb Z}
[N(n,p) \int_p^{p+1}\varphi(g_t\cdot)dt- \sum_{{\scriptstyle k : S_k f\in [p,p+1[}}
\varphi (g_{S_kf}\cdot)]\longrightarrow^{\mathbb P}0
$$
\rm
Let  $\Phi$ denote the function $\Phi(y)=\int^1_0 \varphi (g_ty) dt$. One has $\int_p^{p+1} \varphi (g_t y) dt =\Phi\circ g_p (y)$.

{\bf Proposition} : \it  Under the hypotheses of the theorem, there exists two independent sequences  of i.i.d  centered random variables  with values in $\mathbb Z$ , $(X_i),\ (\xi_i)$ such that :
$$\mathbb E(exp (it\frac{1}{n^{3/4}}  \displaystyle\mathop{\sum}_{p\in \mathbb
Z}N(n,p) \Phi \circ g_p))-\mathbb E(exp (it \frac{1}{n^{3/4}}
\displaystyle\mathop{\sum}_{k=0}^{n-1} \xi_{X_1+\ldots+X_2}))
\displaystyle\mathop{\longrightarrow}_{n\rightarrow +\infty}0. $$

\rm

To prove these two  propositions we use the three following lemmas.

{\bf Lemma 1} : (Moderated deviation)\it 
There exist $C>0, c>0$ such that, for very $n$, every
 $\beta\in ]0,\frac{1}{2}[$, one has
 $$\mathbb P(S_n f>n^{1-\beta})\leq C e^{-cn^{1-2\beta}}$$
\rm
(Here the probability $\mathbb P$ is the Lebesgue measure on
 $\T^2$).
 
{\bf Lemma 2} : (Multiple mixing for the geodesic flow)  \it 
 Let $m$ and $m'$ be two integers,
 $(\Phi_i)_{i=i}^{m+m'}$ h\"older-continuous functions  defined
 on $G/\Gamma, t_1,\leq \ldots\leq t_m\leq 0\leq s_1\ldots\leq s_{m'}, T>0$ real numbers. There
 exists $C>0$, and $\delta\in]0,1[$,such that

 Cov $(\displaystyle\prod_{i=1}^{m}\Phi_i \circ g_{t_i},
 \displaystyle\prod_{j=1}^{m'} \Phi_j\circ g_{s_i+T})\leq
 C(\displaystyle\prod_{i=1}^{m+m'}\| \Phi_i\|_{\infty}+
 \displaystyle\mathop{\sum}_{j}[\Phi_j]\displaystyle\mathop{\sum}
_{i\neq j} \|\phi_j\|_{\infty})\delta^T.$
\rm

{\bf Lemma 3} : \it  Let

- $I$ be an interval of length  1,

- $\varepsilon >0$ a real number,

- $J,K$ be two subintervals of $I$ of length $\frac{1}{[n^{\varepsilon}]}$.

Let $N(n,I)$ denote the quantity :
$$N(n,I)=\sharp \{k\in \{0, -n-1\}/S_k f\in I\}.$$

We define in the same way $N(n,K)$ and $N(n,J)$.

Under the hypotheses of the theorem there exist $\xi>0, C>0$ such that :

$\mathbb E(N(n,I)^2)-[n^{\varepsilon}] \mathbb E (N(n,I)
N(n,J))|\leq C n^{1-\xi}$

$\mathbb E(N(n,J)N(n,K)) [n^{\varepsilon}]^2- \mathbb E (N(n,I))|\leq C n^{1-\xi}$

 $\mathbb E(N(n,I)^3)|\leq C n^{3/2}$
\rm

For the first and second lemmas see \cite{2} et \cite{9}. The third one is deduced from a result of speed of convergence in the local limit theorem for the sums $S_n f$.

\end{document}